\documentclass[ejs,preprint,dvips,twoside,linksfromyear]{imsart}

\RequirePackage[OT1]{fontenc}
\RequirePackage{amsthm,amsmath}
\RequirePackage[numbers]{natbib}
\RequirePackage[colorlinks,citecolor=blue,urlcolor=blue]{hyperref}
\usepackage{srcltx}
\usepackage{graphicx,color,amssymb,latexsym,amsmath,amsfonts}
\usepackage{mathrsfs}

\startlocaldefs
\numberwithin{equation}{section}
\theoremstyle{plain}
\newtheorem{thm}{Theorem}[section]

\newtheorem{lemma}{Lemma}[section]

\newtheorem{corollary}{Corollary}[section]
\endlocaldefs
\begin{document}
\begin{frontmatter}

\title{On Distinguishability of  Hypotheses}

\runtitle{On distinguishability of hypotheses}

\begin{aug}
\author{\fnms{Mikhail} \snm{Ermakov}\ead[label=e1]{erm2512@gmail.com}}
\address{Institute of Problems of Mechanical Engineering, RAS and\\
St. Petersburg State University, St. Petersburg, RUSSIA\\
\printead{e1}}

\thankstext{t2}{Research partially supported by RFFI Grant 11-10-00577.}
\runauthor{M. Ermakov}
\end{aug}

\begin{abstract}
 We consider the problems of hypothesis testing on a probability measures of independent samples, on solutions of ill-posed problems in Gaussian noise, on deconvolution problems  and on   mean measures of Poisson processes. For all these
setups necessary conditions and sufficient conditions are given for distinguishability of sets of hypothesis.
In the case of hypothesis testing on a probability measure and on Poisson mean measure the results are given in terms of weak topology and
topology of weak convergence of measures on all Borel sets.
 The problem of discernibility of hypothesis is also studied. In  other cases the necessary and sufficient conditions of distinguishability are given if the sets of hypotheses are bounded  sets in $L_2$
\end{abstract}



\received{\smonth{8} \syear{2013}}


\end{frontmatter}

\section{\bf  Introduction}The sets of hypothesis and alternatives are distinguishable if there exists uniformly consistent family of tests.  In parametric hypothesis testing   we almost never faced with the problem of distinguishability.
In nonparametric hypotheses testing the problem of distinguishability usually emerges if  we consider the test properties for nonparametric sets of alternatives.

For existence of consistent  nonparametric estimator we usually suppose  certain compactness conditions (see Le Cam and Schwartz \cite{les}, Ibragimov and Khasminskii \cite{ih}, Pfanzagl \cite{pf}, Le Cam \cite{le86}, Yatracos \cite{ya}).
In nonparametric hypothesis testing the  distinguishability of fixed sets of hypotheses  usually does not require such  strong assumptions.
 By this reason on   distinguishability conditions we do not pay  such a serious attention.

The most wellknown papers on distinguishability are the papers of Hoefding, Wolfowitz \cite{ho} and Le Cam, Schwartz \cite{les}.
Hoefding and Wolfowitz  \cite{ho} have proposed a simple form of sufficient conditions of distinguishability. The distinguishability conditions were given in terms of intervals of monotonicity of differences of distribution functions. Le Cam and Schwartz \cite{les} (see also Le Cam \cite{le73},\cite{le86}) have found necessary and sufficient conditions of distinguishability. There are a lot of papers studying special problems related with distinguishability. We  discuss later their relationship with certain results of the paper.

For the problems of hypothesis testing in functional spaces the first result on indistinguishability  has been obtained Le Cam \cite{le73}. For the problem of hypothesis testing
on a density Le Cam \cite{le73} has proved that the center of the ball and the interior of the ball in $L_1$ are indistinguishable. For  signal detection in Gaussian white noise Ibragimov and Hasminski \cite{ih} has shown  that we could not distinguish the hypotheses if the sets of alternatives is the interior of ball in $L_2$. Ingster \cite{in93}
and Ingster, Kutoyants \cite{ik} have obtained similar results for  hypothesis
testing on a density  and for an intensity of Poisson process  respectively. Burnashev \cite{bu} has proved the indistinguishability of
the interior of the ball in $L_p$- spaces, $p> 0$ in the problem of signal detection in Gaussian white noise.
Ermakov  \cite{er00} found necessary and sufficient conditions of distinguishability in $L_2$ for the problem of signal detection in Gaussian white noise and  for  hypothesis testing on a density.

The distinguishability conditions in Le Cam and Schwartz \cite{les} were given  in terms of weak
topology defined on the products of probability measures for all sample sizes. The weak topology was generated all bounded measurable functions.
The goal of this paper is to obtain simple necessary conditions and simple sufficient conditions for distinguishability in terms of weak topologies defined on probability measures. We consider only the case of fixed sets of hypotheses. The approaching sets of hypotheses  should satisfy the same indistinguishability conditions.

In $L_2$-spaces, for all setups, the necessary and sufficient conditions of distinguishability are given.
For more general setups of hypothesis testing on a probability measure and on a Poisson mean measure
we can find necessary and sufficient conditions of distinguishability only
for relatively compact sets of measures in the topology of setwise convergence on
all Borel sets ($\tau$-topology).
These theorems are based on wellknown results on compactness in this topology (see Gaenssler \cite{ga}).
The distinguishability conditions for these setups are given for more general  assumptions. All distinguishability results are rather evident.

The setup of the paper is the following.  Let we be given a family of statistical experiments $  \frak{E}_\epsilon = (\Omega_\epsilon,  \frak{B}_\epsilon,  \frak{P}_\epsilon)$ where $(\Omega_\epsilon,  \frak{B}_\epsilon)$  is sample space with $\sigma$-fields of Borel sets $ \frak{B}_\epsilon$ and  $ \frak{P}_\epsilon= \{P_{\theta,\epsilon}, \theta \in \Theta\}$  is a family of probability measures.
One needs to test a hypothesis $H_0: \theta \in \Theta_0 \subset \Theta$ versus alternative $H_1: \theta \in \Theta_1 \subset \Theta$. In what follows, $\epsilon \in R^1, \epsilon > 0$ or $ \epsilon = n^{-1/2}$ with $n =1,2,\ldots$.

For any family of tests $K_\epsilon$ denote $\alpha_{\theta}(K_\epsilon), \theta \in \Theta_0$, and $\beta_{\theta}(K_\epsilon), \theta \in \Theta_1,$ their type I and type II error probabilities respectively.

Denote
$$
\alpha(K_\epsilon) = \sup_{\theta \in \Theta_0} \alpha_{\theta}(K_\epsilon) \quad
\mbox{\rm and}
\quad
\beta(K_\epsilon) = \sup_{\theta \in \Theta_1} \beta_{\theta}(K_\epsilon).
$$
A family of tests $K_\epsilon$, is called consistent (see Lehmann \cite{le}, van der Vaart \cite{van}
and references therein), if
$$
\limsup_{\epsilon\to 0}\alpha_{\theta_0}(K_\epsilon)  = 0\quad
\mbox{\rm and}
\quad \limsup_{\epsilon\to 0}\beta_{\theta_1}(K_\epsilon)  = 0
$$
for all $\theta_0\in \Theta_0$ and $\theta_1 \in \Theta_1$.

We are interested in the existence of uniform consistent tests (see Hoefding and Wolfowitz \cite{ho}).
We say that a family of tests $K_\epsilon, \alpha(K_\epsilon) < \alpha$  is uniformly consistent if
$$
\lim_{\epsilon\to 0} \beta(K_\epsilon) = 0
$$
for all $0 < \alpha <1$.

We say that   hypothesis $H_0$ and  alternative $H_1$ are distinguishable (see Hoefding and Wolfowitz \cite{ho}) if  there is  uniformly  consistent family of tests.

 We say that hypotheses $H_0$ and  alternative $H_1$ are indistinguishable (see Hoefding and Wolfowitz \cite{ho}) if for each $\epsilon$ for each family of tests $K_\epsilon$ there holds
$$\alpha(K_\epsilon) + \beta(K_\epsilon) \ge 1$$.

The paper is organized as follows. Section \ref{s2} contains  main results. In subsection \ref{s2.1} necessary conditions and sufficient conditions for distinguishability of sets of probability measures of independent sample are given. The problem of discernibility is discussed as well. In subsection \ref{s2.2} we remind the conditions of distinguishability in $L_2$ for the problems of signal detection in Gaussian white noise and hypothesis testing on a density of independent sample (see Ermakov \cite{er00}).   The results of subsection \ref{s2.2} are implemented in subsections \ref{s2.4} and \ref{s2.5}. In subsection \ref{s2.3} we show that similar results of subsections of \ref{s2.1} and \ref{s2.2} are valid for the problem of hypothesis testing on Poisson mean measure. In subsection \ref{s2.4} we study  the distinguishability conditions for ill-posed problems with Gaussian  random noise and the problem of signal detection in the heteroscedastic Gaussian white noise.  In subsection \ref{s2.5} we find the distinguishability conditions in deconvolution problem. The proofs of results of subsections \ref{s2.1} and \ref{s2.3} are given in sections \ref{s3} and \ref{s4} respectively.

We  shall denote by letters $c$ and $C$ generic constants. Denote $1_A(x)$ the indicator of the set $A$. Denote $[a]$ the whole part of $a \in R^1$. For any measures $P_1, P_2$ denote $P_1\otimes P_2$ the product measure.
\section{\bf Main Results \label{s2}}
\subsection{Hypothesis testing on a probability measure of independent sample \label{s2.1} } Let $X_1,\ldots,X_n$ be i.i.d.r.v.'s on a probability space $(\Omega, \frak{B}, P)$
 where $\frak{B}$ is $\sigma$-field of Borel sets on Hausdorff topological space $\Omega$.
 Denote $\Theta$ the set of all probability measures on $(\Omega, B)$.

The coarsest topology in $\Theta$ providing the continuous mapping
 $$
 P \to P(B), \qquad P \in \Theta
 $$
 for all $B \in \frak{B}$ is called the $\tau$-topology (see Groeneboom,  Oosterhoff, Ruymgaart \cite{gor}  and   Dembo,   Zeitouni \cite{dem}) or the topology of setwise convergence on all Borel sets
  (see Ganssler \cite{ga} and Bogachev \cite{bo}).

 For any set $\Psi \in \Theta$ denote $\frak{cl}_\tau(\Psi)$ the closure of $\Psi$ in the $\tau$-topology.

 We say that the set $\Psi \subset \Theta$ is relatively compact in $\tau$-topology if the closure of $\Psi$ is compact in $\tau$-topology.
 \begin{thm}\label{t3.1} There hold.

   {\sl i.} Let  $\Theta_0$ and $\Theta_1$ be relatively compact in $\tau$-topology. Then hypothesis $H_0$ and  alternative $H_1$ are indistinguishable if $\frak{cl}_\tau(\Theta_0) \cap \frak{cl}_\tau (\Theta_1) \ne \emptyset$.

   {\sl ii.} If $\Theta_0$ or $\Theta_1$ is relatively compact in $\tau$-topology and $\frak{cl}_\tau(\Theta_0) \cap \frak{cl}_\tau (\Theta_1) = \emptyset$, then the hypothesis $H_0$ and alternative $H_1$ are distinguishable. There exists a sequence of tests $K_n$
 and constant $n_0$ such that, for all $n> n_0$, we have
 \begin{equation}\label{m1}
 \alpha(K_n) \le \exp \{-cn\} \quad \mbox{and} \quad \beta(K_n) \le \exp\{-cn\}.
\end{equation}
 \end{thm}
 Le Cam \cite{le73} and Schwartz \cite{sch} have shown that, if the hypothesis and alternative are not indistinguishable,
  then (\ref{m1}) holds. Thus, if the hypothesis and alternative are not indistinguishable,  they are distinguishable. This statement has been proved Hoefding and Wolfowitz \cite{ho}. The reasoning in  Le Cam \cite{le73} and Schwartz \cite{sch} easily can be  extended on all models of hypothesis testing in  this paper.

  It should be mentioned that the property of exponential decay of type I and type II error probabilities was studied
  in a large number of papers (see Hoefding and Wolfowitz \cite{ho}, Schwartz \cite{sch},
 Le Cam \cite{le73}, Dembo, Zeitouni \cite{dem}, Barron \cite{ba}, Ermakov \cite{er93} and references therein).
 In the proof of Theorem \ref{t3.9} we implement (\ref{m1}) for a sequence of tests defined explicitly. This sequence of tests
 is defined in the proof of  (\ref{m1}).

 If  set $\Psi$ is relatively compact in $\tau$-topology, then, by Theorem 2.6 in Gaenssler \cite{ga}, the set $\Psi$ is equicontinuous and
 there exists probability measure $\nu$ such that $P<<\nu$ for all $P \in \Psi$.
  This implies that for any $\delta>0$ there exists $\epsilon>0$ such that, if $\nu(B)<\epsilon, B\in \frak{B}$, then $P(B)<\delta$ for all
  $P \in \Psi$. The set of densities $\frak{F} = \left\{f: f = \frac{dP}{d\nu}, P \in \Psi \right\}$ is uniformly integrable.

  If the sequence of probability measures $P_m$ converges to $P\in\Theta$ in $\tau$-topology, then this
  sequence is relatively compact in $\tau$-topology (see Ganssler \cite{ga}). Thus the condition of indistinguishability
   can be given in the following form.

   For any set $\Psi \in \Theta$ denote $\mbox{cl}_{s\tau}(\Psi)$ the sequential closure of $\Psi$ in $\tau$-topology.
   Then, for any sets $\Theta_0, \Theta_1 \subset \Theta$,  $\mbox{cl}_{s\tau}(\Theta_0) \cap \mbox{cl}_{s\tau} (\Theta_1) \ne \emptyset$ implies
    indistinguishability of $H_0$ and $H_1$. We could not prove indistinguishability for any sets $\Theta_0, \Theta_1 \subset \Theta$
   if $\mbox{cl}_\tau(\Theta_0) \cap \mbox{cl}_\tau (\Theta_1) \ne \emptyset$.
   The map $P \to P\times P$ of $\Theta \to \Theta\times\Theta$ is not continuous in the $\tau$-topology (see 8.10.116 in Bogachev \cite{bo}).

 {\sl Example} 2.1. Let $\nu$ is  Lebesgue measure in $(0,1)$ and let we consider the problem of hypothesis  testing on a density $f$ of probability measure $P$. Let $H_0 : f(x) =1,\, x \in (0,1)$ and $\Theta_1 = \{f_1, f_2, \ldots\}$ with $f_i(x) = 1 + \sin(2\pi i x), \, x \in (0,1)$.

 For any measurable set $B \in \frak{B}$ we have
  $$
 \lim_{i\to\infty} \int_B f_i(x) \, dx = \int_B dx.
 $$
 Therefore $H_0$ and $H_1$ are indistinguishable.

 For any probability measures $P<<\nu$ and $Q<<\nu$ define the variational distance
$$
\mbox{\rm var}(P,Q) = \frac{1}{2}\int_\Omega |dP/d\nu - dQ/d\nu| \,d\nu
$$
For any sets of probability measures $A$ and $B$ denote
$$
\mbox{\rm var}(A,B) = \inf\{\mbox{var}(P,Q):\, P\in A, Q \in B\}.
$$
Denote $[A]$  the convex hull of set $A \subset\Theta$.

Kraft \cite{kr} proved the following Theorem.
\begin{thm}\label{t5} Let the probability measures in $\Theta_0 \cup \Theta_1$ be absolutely continuous with respect to measure $\nu$. Then, for any test $K$ there holds
\begin{equation}\label{m2}
\alpha(K,\Theta_0) + \beta(K,\Theta_1) \ge 1 - \mbox{\rm var}([\Theta_0],[\Theta_1]).
\end{equation}
\end{thm}
The $L_1$-norm is standard tool in hypothesis testing (see Hoefding and Wolfowitz \cite{ho}, Le Cam \cite{le73, le86}, Lehmann \cite{le},
van der Vaart \cite{van}, Devroye and Lugosi \cite{de}, Ingster and Suslina \cite{is} and references therein).

A sequence of densities $f_k$ converges to density $f_0$ in $L_1(\nu)$  (see Dunford and Schwarz \cite{du}, Th.12 sec.8 Ch IV, or Iosida, \cite{io} Th.5  sec.1 Ch. V)
if for any measurable set $B \in B$ there hold
 \begin{equation}\label{m3}
 \lim_{k\to\infty} \int_B f_k\, d\nu= \int_B f_0 d\nu.
\end{equation}
 and $f_k$ converges to $f_0$ in measure.

 If (\ref{m3}) holds and $f_k$ does not converges to $f_0$ in measure we could not distinguish
 the set of hypotheses $\Theta_0=\{f_0\}$ and the set of alternatives $\{f_1,f_2,\ldots\}$.
 By Mazur Theorem (see Iosida \cite{io}, Th.2 sec.1 Ch.5 ), weak convergence $f_k$ to $f_0$ implies convergence of convex combinations of $f_k$ to $f_0$ in $L_1(\nu)$. Therefore the right-hand side of (\ref{m2}) equals one.

 The proof of distinguishability in Theorem \ref{t3.1} is based on the statement of Theorem \ref{t3.3}.
 \begin{thm}\label{t3.3} Let $\Theta_0$ or $\Theta_1$ be relatively compact in the $\tau$-topology and $\frak{cl}_\tau(\Theta_0) \cap \frak{cl}_\tau (\Theta_1) = \emptyset$. Then there are measurable sets $B_1,\ldots,B_k \in \frak{B}$ such that the closures of the sets $\Theta_{0k} =\{u: u = (P(B_1),\ldots,P(B_k)), P \in \Theta_0\}$ and $\Theta_{1k} =\{v: v = (Q(B_1),\ldots,Q(B_k)), Q \in \Theta_1\}$  have no common points.
 \end{thm}
 \vskip 0.3cm
 {\sl Remark 2.1.} Let $\alpha(K) + \beta(K) < 1-\delta, \delta>0,$ for  test $K(X_1,\ldots,X_m)$.
 For any $\epsilon > 0, \epsilon < \delta/2$ there is simple function
 $$
 L_\epsilon(x) = \sum_{i=1}^k a_i 1_{B_i}(x), \quad a_i \in [0,1],\,\, B_i \in \frak{B}^m,
 \,\,x  \in \Omega^m = \Omega\times\ldots\times\Omega
 $$
 such that
 \begin{equation}\label{lu1}
 \sup_{x\in\Omega^m} |L_\epsilon(x)- K(x)| < \epsilon.
 \end{equation}
 Then $\alpha(L) +\beta(L) < 1 - \delta + 2\epsilon < 1$.

 Thus, if the sets of hypotheses and alternatives are distinguishable, Theorem \ref{t3.3} should be valid
 at least for the Borel sets $B_1,\ldots,B_k$ in $\Omega^m$ for some sample size $m$.
\vskip 0.3cm

 In probability and statistics the traditional technique is  weak topology generated continuous functions. Below   versions of Theorems \ref{t3.1} and \ref{t3.3} are given in terms of weak topology.

Let $\Omega$ be Polish space.
Let $\bold{C}$ be the set of all real bounded continuous functions $f: \Omega \to R^1$. We say that $P_k \in \Theta $  converges to $P \in \Theta$ in weak topology if
$$
\lim_{k\to \infty} \int_\Omega f \, dP_k = \int_\Omega f \, dP
$$
for  each $f \in \bold{C}$.

For any set $\Psi \subset \Theta$ we denote $\frak{cl}_c(\Psi)$ the closure of $\Psi$ in the weak topology.
 \begin{thm}\label{t3.5} There hold.

 {\sl i.} Let $\Theta_0$ and $\Theta_1$ are relatively compact in $\tau$-topology. Then the hypothesis $H_0$ and the alternative $H_1$ are indistinguishable if $\frak{cl}_c(\Theta_0) \cap \frak{cl}_c (\Theta_1) \ne \emptyset$.

 {\sl ii.} If $\Theta_0$ or $\Theta_1$ is relatively compact in weak topology and $\frak{cl}_c(\Theta_0) \cap \frak{cl}_c (\Theta_1) = \emptyset$ then the hypothesis $H_0$ and alternative $H_1$ are distinguishable. There exists a sequence of tests $K_n$
 and constant $n_0$ such that, for all $n> n_0$, we have
 \begin{equation}\label{m4}
 \alpha(K_n) \le \exp \{-cn\} \quad \mbox{and} \quad \beta(K_n) \le \exp\{-cn\}.
\end{equation}
 \end{thm}
 The proof of distinguishability in Theorem \ref{t3.5} is based on the following Theorem \ref{t3.7}.
  \begin{thm}\label{t3.7}  Let $\Theta_0$ or $\Theta_1$ is relatively compact in weak topology and $\frak{cl}_c(\Theta_0) \cap \frak{cl}_c (\Theta_1) = \emptyset$. Then there exist open sets $B_1,\ldots,B_k \subset\Omega$ such that the closures of the sets $\Theta_{0k} =\{u: u = (P(B_1),\ldots,P(B_k)), P \in \Theta_0\}$ and $\Theta_{1k} =\{v: v = (Q(B_1),\ldots,Q(B_k)), Q \in \Theta_1\}$  have no common points.
 \end{thm}

 Theorems \ref{t3.1},\ref{t3.3} and \ref{t3.5},\ref{t3.7} can easily be extended on the  problems of hypothesis testing with several independent samples.

 Let $X_1,\ldots,X_n$ and $Y_1,\ldots,Y_n$ be i.i.d.r.v.'s having measures $P_1$ and $P_2$ respectively. The probability measures $P_1$ and $P_2$ are defined on   $\sigma$-field $\frak{B}$ of Borel sets in Hausdorff topological space $\Omega$. The problem is to test the hypothesis $H_0: P = P_1\otimes P_2 \in \Theta_0 \subset \Theta\otimes \Theta$ versus the alternative $H_1: P = P_1\otimes P_2 \in \Theta_1 \subset \Theta\otimes \Theta$ where $\Theta_0\cap\Theta_1 = \emptyset$.

 For this setup we can define $\tau$-topology and weak topology
 as the corresponding product topologies. After that Theorems \ref{t3.1},\ref{t3.3} and \ref{t3.5},\ref{t3.7} become valid for the problems of hypothesis testing with several independent samples in terms of these product topologies. The proofs of these statements are similar.

  Devroye, Lugosi \cite{de} and Dembo, Peres \cite{dep} studied the problem of discernibility of hypotheses.
 The sets of hypotheses $H_0: P \in \Theta_0$ and alternatives $H_1: \theta \in \Theta_1$  are called discernible if there is a sequence of tests $K_n$ such that
 \begin{equation}\label{m5}
 P(K_n = 1 \quad\mbox{for only finitely many}\quad n) =1 \quad\mbox{for all}\quad P \in \Theta_0
 \end{equation}
  and
 \begin{equation}\label{m6}
 P(K_n = 0\quad \mbox{for only finitely many}\quad n) =1 \quad\mbox{for all}\quad P \in \Theta_1.
 \end{equation}
  Kulkarni and Zeitouni \cite{ku} considered a slightly different setup.

 It is  clear that (\ref{m5}) and (\ref{m6}) implies the consistency of tests $K_n$. However (\ref{m5}) and (\ref{m6}) does not imply uniform consistency of $K_n$. The discernibility does not imply the distinguishability. The hypothesis and alternatives in Example 2.1 are discernible (see Theorem 2 in  Dembo and Peres \cite{dep}). Below we consider the problem of uniform discernibility.

 We say that positive random variable $N$ is  uniformly bounded  if
 \begin{equation}\label{m9}
 \lim_{C\to\infty}\sup_{P\in\Theta_0\cap\Theta_1} P(N > C) = 0
 \end{equation}
 We say that the sets of hypotheses $H_0: P \in \Theta_0$ and alternatives $H_1: \theta \in \Theta_1$  are uniformly discernible if there are sequence of tests $K_n$ and uniformly bounded positive random variable $N$ such that
 \begin{equation}\label{m7}
 P(K_n = 0\quad \mbox{for all} \quad n>N) =1 \quad\mbox{for all}\quad P \in \Theta_0
 \end{equation}
 and
 \begin{equation}\label{m8}
 P(K_n = 1\quad \mbox{for all} \quad n>N) =1 \quad\mbox{for all}\quad P \in \Theta_1.
 \end{equation}
  We  say that the sequence of tests $K_n$ uniformly discern the hypothesis $H_0$ and alternative $H_1$.

 It is clear that uniform discernibility implies distinguishability.
\begin{thm}\label{t3.9} Let $\Omega$ be Hausdorff topological space and
let $\frak{B}$ be the set of Borel sets in $\Omega$.
Let  hypothesis $H_0 : P \in \Theta_0$ and  alternative $H_1: \theta \in \Theta_1$
are distinguishable. Then the hypothesis $H_0$ and the alternative $H_1$  are uniformly discernible
and there exists  sequence of tests $K_n$ and $t > 0$ such that
\begin{equation}\label{m10}
\sup_{P\in\Theta_0\cap\Theta_1} E_P [\exp\{tN\}] < C < \infty.
 \end{equation}
\end{thm}
Note that exponential decay of type I and type II error probabilities (see (\ref{m1}) and (\ref{m4}))
follows from (\ref{m10}).

The hypothesis $H_0: P\in \Theta_0$ is often considered for alternatives $H_1: P \in \Theta_1= \Theta\setminus\Theta_0$.

We say that the hypothesis $H_0$ are uniformly asymptotically discernible, if there are a sequence of sets $\Theta_k, \Theta_k \subseteq \Theta_{k+1}$,  $\cup_{k=1}^\infty\Theta_k \cup \frak{cl}_c(\Theta_0) = \Theta$,  and sequence of tests $K_n=K_n(X_1,\ldots,X_n)$  such that, for any $k, k =1,2,\ldots$, the sequence of tests $K_n$ uniformly discern the hypothesis $H_0$ and alternative $H_k: P \in \Theta_k$. In this case we call the set $\Theta_0$ uniformly asymptotically discernible.
\begin{thm}\label{t3.10} Let $\Omega$ be Polish space. Let $\Theta_0$ be relatively compact in weak topology. Then $\Theta_0$ is uniformly asymptotically discernible.
 \end{thm}
\subsection{Signal detection and hypothesis testing on a density in $L_2$ \label{s2.2}}
The problem of signal detection is treated  in the following setup.
Let we observe a realization of stochastic process $Y_\epsilon(t), t \in (0,1)$ defined by the stochastic differential equation
$$
dY_\epsilon(t) = S(t) dt + \epsilon dw(t)
$$
where $S$ is unknown signal and $dw(t)$ is Gaussian white noise.

The problem of hypothesis testing on a density $S  \in \Theta\subset L_2(\nu)$ is studied in the following setup. Let $X_1,\ldots,X_n$ be i.i.d.r.v.'s defined on a probability space $(\Omega, \frak{B}, P)$.
Suppose the probability measure $P$ is continuous w.r.t. probability measure $\nu$ and has the density $S(x)= dP/d\nu (x), x \in \Omega$.

The problem is  to test the hypothesis on unknown parameter $\theta = S$ with $\Theta_0, \Theta_1 \subset \Theta$
where $\Theta$ is closed subspace of $L_2(\nu)$.

For any subspace $\Gamma \subset L_2(\nu)$ denote  $\Pi_\Gamma$ the projection operator on the subspace $\Gamma$. For any set $\Psi\subset L_2(\nu)$ denote $\bar\Psi$ the closure of $\Psi$ in $L_2(\nu)$.

\begin{thm} \label{t1}  Suppose the sets $\Theta_0$ and $\Theta_1$ are  bounded  in $L_2(\nu)$.
Then both for the problems of testing hypotheses on a density and signal detection the hypothesis $H_0$
and alternative $H_1$ are

{\sl i.} distinguishable if there exists a finite dimensional subspace $\Gamma \subset \Theta$
such that $\Pi_\Gamma\bar\Theta_0 \cap \Pi_\Gamma\bar\Theta_1 = \emptyset$

{\sl ii.} indistinguishable if there does not exist a finite dimensional subspace $\Gamma \subset \Theta$
such that $\Pi_\Gamma\bar\Theta_0 \cap \Pi_\Gamma\bar\Theta_1 = \emptyset$.
\end{thm}
The proof of Theorem \ref{t1} was given in Ermakov \cite{er00}. Similar proof of Theorem \ref{t4.3} on the problem of hypothesis on intensity of Poisson process  is given in this paper with simplified reasoning.
\subsection{Hypothesis testing on a  mean measure of Poisson random process \label{s2.3}} Let we be given $n$ independent realizations $\kappa_1,\ldots,\kappa_n$ of Poisson random process with mean measure $\mu$  defined on Borel sets $\frak{B}$ of Hausdorff space $\Omega$.
The problem is to test a hypothesis $H_0: \mu \in \Theta_0\subset\Theta$ versus $H_1: \mu \in \Theta_1\subset\Theta$ where  $\Theta$ is the set of all measures $\mu, \mu(\Omega) < \infty$.

The results for this setup are the same as in the problem of hypothesis testing on probability measure of independent sample  (compare Theorem \ref{t3.1} and  Theorem \ref{t4.1}).
\begin{thm}\label{t4.1} There hold.

   {\sl i.} Suppose the sets $\Theta_0$ and $\Theta_1$ are relatively compact in $\tau$-topology.
   Then hypothesis $H_0$ and alternative $H_1$ are indistinguishable if $\frak{cl}_\tau(\Theta_0) \cap \frak{cl}_\tau(\Theta_1) \ne \emptyset$.

   {\sl ii.} If $\Theta_0$ or $\Theta_1$ is relatively compact in $\tau$-topology and $\frak{cl}_\tau(\Theta_0) \cap \frak{cl}_\tau (\Theta_1) = \emptyset$, then the hypothesis $H_0$ and alternative $H_1$ are distinguishable. There exists a sequence of tests $K_n$
 and constant $n_0$ such that, for all $n> n_0$, we have
 \begin{equation}\label{pu5}
 \alpha(K_n) \le \exp \{-cn\} \quad \mbox{and} \quad \beta(K_n) \le \exp\{-cn\}.
\end{equation}
 \end{thm}
 The statement of Theorem \ref{t3.3} is valid  for the setup of Theorem \ref{t4.1} as well.
\begin{thm}\label{t4.2} For the problem of hypothesis testing on Poisson mean measure the statements of Theorems \ref{t3.5}, \ref{t3.7} and \ref{t3.9}, \ref{t3.10} are valid with the sets $\Theta_0, \Theta_1$ belonging  to the set  $\Theta$  of all measures $\mu, \mu(\Omega) < \infty$.
\end{thm}
Let $\nu$ be Poisson mean measure and let $\Theta$ be a subset  of the set of Poisson intensity functions $\lambda = d\mu/d\nu, \mu << \nu$. Suppose that $\Theta$ is closed
subspace of $L_2(\nu)$ and $\nu(\Omega) < \infty$.  Suppose $\lambda_c(t) \equiv 1 \in \Theta$. For this setup the statement of Theorem \ref{t1} holds.
\begin{thm} \label{t4.3}  Let the sets $\Theta_0$ and $\Theta_1$ be  bounded in $L_2(\nu)$. 
Then,  for the problem of testing hypotheses on intensity of Poisson random process, the sets of hypotheses $\Theta_0$
and alternatives $\Theta_1$ are

 {\sl i.} indistinguishable if there does not exist a finite dimensional subspace $\Gamma \subset \Theta$
 such that $\Pi_\Gamma\bar\Theta_0 \cap \Pi_\Gamma\bar\Theta_1 = \emptyset$

 {\sl ii.} distinguishable if there  exists a finite dimensional subspace $\Gamma \subset \Theta$
 such that $\Pi_\Gamma\bar\Theta_0 \cap \Pi_\Gamma\bar\Theta_1 = \emptyset$.
\end{thm}
 \subsection{Hypothesis testing on a solution of ill-posed problem \label{s2.4}} In Hilbert space $H$
  we wish to test a hypothesis on a  vector $\theta = x \in \Theta \subset H$ from
  the observed  Gaussian random vector
  $$
  Y = Ax + \epsilon \xi.
  $$
  Hereafter $A: H \to H$ is known linear operator and $\xi$ is Gaussian random vector having known covariance operator $R: H \to H$ and
  $E\xi = 0$.

  Suppose that the kernels of $A$ and $R$ equal zero.

  Denote $\Pi_\Gamma$ the projection operator on linear subspace $\Gamma\subset H$.
  For any set $\Psi \subset H$ denote $\bar\Psi$ the closure of $\Psi$ in $H$.

  The set $\Theta$ is closed linear subspace of $H$ and $\Theta_0, \Theta_1 \subset \Theta$.

  Define the sets $\Lambda_0 = R^{-1/2}A\Theta_0$ and $\Lambda_1=R^{-1/2}A\Theta_1$.  Denote $\Lambda$
 the closure of linear space generated
 linear combinations of vectors from $\Lambda_0\cup\Lambda_0$.
\begin{thm}\label{t21} Let the sets $\Lambda_0$ and $\Lambda_1$ are bounded in $H$. Then the hypothesis $H_0$ and the alternative $H_1$ are

{\sl i.} indistinguishable if there does not exist finite dimensional subspace $\Gamma \subset \Lambda$ such
that $\Pi_\Gamma\bar\Lambda_0 \cap \Pi_\Gamma\bar\Lambda_1=\emptyset$

{\sl ii.} distinguishable if there exists finite dimensional subspace $\Gamma \subset \Lambda$ such
that $\Pi_\Gamma\bar\Lambda_0 \cap \Pi_\Gamma\bar\Lambda_1=\emptyset$.
\end{thm}
{\it Proof of Theorem \ref{t21}.} Denote $Z=  R^{-1/2}Y$ and $z = EZ$. Consider the problem of testing hypothesis $\bar H_0: z \in \Lambda_0$ versus $\bar H_1: z \in \Lambda_1$.
The random vector $Z - EZ$ is Gaussian white noise.
Therefore, by Theorem \ref{t1}, the hypothesis $\bar H_0$ and the alternative  $\bar H_1$ are distinguishable  iff there exists finite dimensional subspace
$\Gamma\subset\Lambda$ such that $\Pi_\Gamma\bar\Lambda_0 \cap \Pi_\Gamma\bar\Lambda_1 = \emptyset$  that completes the proof of Theorem \ref{t21}.
\begin{thm}\label{t22}  Let the sets $\Lambda_0$ and $\Lambda_1$ are bounded in $H$. Then the hypothesis $H_0$ and the alternative $H_1$ are

{\sl i.} indistinguishable if there does not exist a finite dimensional subspace $\Psi \subset \Theta$
such that $\Pi_\Gamma\bar\Theta_0 \cap \Pi_\Psi\bar\Theta_1 = \emptyset$

{\sl ii.} distinguishable if there exists a finite dimensional subspace $\Psi \subset \Theta$
such that $\Pi_\Psi\bar\Theta_0 \cap \Pi_\Psi\bar\Theta_1 = \emptyset$.
\end{thm}
{\it  Proof of Theorem \ref{t22}.} Let the hypothesis $H_0$ and alternative $H_1$ be distinguishable. Let the space $\Gamma$ be the same as in the proof of Theorem \ref{t21}. The operator $\Pi_\Gamma R^{-1/2}A$ has values in finite dimensional linear subspace $\Gamma$. Denote $U = \{x: \Pi_\Gamma R^{-1/2}Ax = 0, x \in H\}$ the kernel of operator $\Pi_\Gamma R^{-1/2}A$. Denote $\Psi = \{y: (x,y) = 0, x\in U, y \in H\}$. The linear subspace $\Psi$ is finite dimensional. Since $\Pi_\Gamma\bar\Lambda_0 \cap \Pi_\Gamma\bar\Lambda_1=\emptyset$ then $\Pi_\Gamma\bar\Theta_0 \cap \Pi_\Psi\bar\Theta_1 = \emptyset$. This completes the proof of Theorem \ref{t22}.
\vskip 0.25cm
\noindent{\sl Remark.} Suppose
the operator $R^{-1/2}A$ be bounded. Suppose the sets $\Theta_0$ and $\Theta_1$ are bounded  in $H$ Suppose the set $\Theta_0$ of hypotheses and the set $\Theta_1$ of alternatives are distinguishable.  Define the operator  $D= R^{-1/2}A$.

For any operator $G: H\to H$ define the norm of the operator
$$
||G|| = \sup_{||x||=1} |(x,Gx)|.
$$
Hereafter  $||x||$ denotes the norm of $x \in H$.

For any $\delta > 0$ we can represent the operator $D$ as  sum of two self-adjoint operators
\begin{equation*}
D = D_{1\delta} + D_{2\delta}
\end{equation*}
where
$$
|| D_{1\delta}^{-1}|| = \delta^{-1} \quad \mbox{and} \quad  ||D_{2\delta}|| = \delta.
$$
Denote $U_\delta$ the kernel of operator $D_{1\delta}$ and $V_\delta = \{y: (x,y) = 0, x\in U_\delta, y \in H\}$.

By Theorem \ref{t21} there exists a finite dimensional subspace $\Upsilon\subset H$ such that $\Pi_\Upsilon\bar\Lambda_0 \cap
\Pi_\Upsilon\bar\Lambda_1 = \emptyset$. Denote $\Phi_0 = \Pi_\Upsilon\bar\Lambda_0$ and $\Phi_1 =\Pi_\Upsilon\bar\Lambda_1$.

For any set $B$ and any $\tau > 0$ denote
\begin{equation*}
B_\tau = \{x: ||x- y|| \le \tau, y \in B, x \in H\}
\end{equation*}
Since the sets $\Phi_0$ and  $\Phi_1$ are closed, there exists $\tau>0$ such that $\Phi_{0\tau}\cap \Phi_{1\tau} = \emptyset$.

Since $\Theta_0$ and $\Theta_1$ are bounded there exists $\delta >0$ such that, for any $x \in \Theta_0$ and any $y \in \Theta_1$ there hold $||D_{2\delta}x|| < \tau$ and $||D_{2\delta}y|| < \tau$.

Therefore $ D_{1\delta}\Theta_0 \subset \Lambda_{0\tau}$ and $ D_{1\delta}\Theta_1 \subset \Lambda_{1\tau}$ and $\Pi_\Upsilon D_{1\delta}\bar\Theta_0 \cap \Pi_\Upsilon D_{1\delta}\bar\Theta_1 = \emptyset$.

Therefore, if we consider the problems of distinguishability of $\Theta_0, \Theta_1$ or  $\Pi_{V_\delta}\Theta_0, \Pi_{V_\delta}\Theta_1$ the result will be the same.
\vskip 0.25cm

The problem of distinguishability for hypothesis testing on a signal  in the heteroscedastic Gaussian white noise can be considered as a particular case of the setups of Theorems \ref{t21} and \ref{t22}.

Let we observe a random process $Y(t), t \in (0,1)$ defined the stochastic differential equation
$$
dY(t) = S(t) dt + \epsilon h(t) dw(t)
$$
where $S$ is unknown signal, $h(t)$ is a weight function and $dw(t)$ is Gaussian white noise.

One needs  to test a hypothesis on unknown parameter $\theta = S$ with $\Theta_0, \Theta_1 \subset \Theta$
where $\Theta$  is closed linear subspace of $ L_2(0,1)$.
\begin{thm}\label{t3} Let $0< c <  h(t) < C < \infty$ for all $t \in (0,1)$. Let  sets  $\Theta_0$ and $\Theta_1$ be bounded  in $L_2(0,1)$. Then the hypothesis $H_0$ and alternative $H_1$ are

 {\sl i.} indistinguishable if there does not exist a finite dimensional subspace $\Gamma \subset \Theta$
 such that $\Pi_\Gamma\bar\Theta_0 \cap \Pi_\Gamma\bar\Theta_1 = \emptyset$

 {\sl ii.} distinguishable if there  exists a finite dimensional subspace $\Gamma \subset \Theta$
 such that $\Pi_\Gamma\bar\Theta_0 \cap \Pi_\Gamma\bar\Theta_1 = \emptyset$.
\end{thm}
In the setup of Theorem \ref{t3} the operators $R$ and $R^{-1}$ induced the weight function $h$ are bounded. By Theorem \ref{t22}, this implies Theorem \ref{t3}.
\subsection{Hypothesis testing on a solution of deconvolution problem \label{s2.5}}
Let we observe i.i.d.r.v.'s $Z_1,\ldots,Z_n$ having  density $h(z), z \in R^1$ with respect to Lebesgue measure. It is known that $Z_i = X_i + Y_i, 1 \le i \le n$ where $X_1,\ldots,X_n$ and $Y_1,\ldots,Y_n$ are i.i.d.r.v.'s with  densities $f(x), x \in R^1$ and $g(y), y \in R^1$ respectively.  The density $g$ is known.

The problem is to test the hypothesis $H_0: f \in \Theta_0$ versus the alternative $H_1: f \in \Theta_1$ where $\Theta_0, \Theta_1 \subset \Theta$. The set $\Theta$ is closed subspace of $L_2(R^1)$.

Suppose
 $g \in L_2(R^1)$.

 Denote
 $$
 (g \ast f) (x) = \int_{-\infty}^{\infty} g(x-y) f(y)
\, dy, \quad x \in R^1
$$
and
$$
\hat g(\omega) = \int_{-\infty}^{\infty} \exp\{i\omega x\} g(y) \, dy, \quad \omega \in R^1.
$$
Define the sets
$\Lambda_i = \{h: h =g\ast f, f \in \Theta_i\}$ with $i = 0,1$ and $\Lambda = \{h: h =g\ast f, f \in \Theta\}$.
\begin{thm}\label{t23} Suppose the sets $\Lambda_0$ and $\Lambda_1$ are bounded  in $L_2(R^1)$. Let $\hat g(\omega) $ be continuous and  let $\hat g(\omega) \ne 0$ for all $\omega \in R^1$. Then the hypothesis $H_0$ and the alternative $H_1$ are

{\sl i.} indistinguishable if there does not exist finite dimensional subspace $\Gamma\subset\Lambda$ such
that $\Pi_\Gamma\bar\Lambda_0 \cap \Pi_\Gamma\bar\Lambda_1=\emptyset$

{\sl ii.} distinguishable if there exists finite dimensional subspace $\Gamma\subset\Lambda$ such
that $\Pi_\Gamma\bar\Lambda_0 \cap \Pi_\Gamma\bar\Lambda_1=\emptyset$.
\end{thm}
\begin{thm}\label{t24} Suppose the sets $\Theta_0$ and $\Theta_1$ are  bounded  in $L_2(R^1)$. Let $\hat g(\omega) $ be continuous and  let $\hat g(\omega) \ne 0$ for all $\omega \in R^1$. Then the hypothesis $H_0$ and the alternative $H_1$ are

{\sl i.} indistinguishable if there does not exist a finite dimensional subspace $\Gamma \subset \Theta$
 such that $\Pi_\Gamma\Theta_0 \cap \Pi_\Gamma\Theta_1 = \emptyset$

{\sl ii.} distinguishable if there exists a finite dimensional subspace $\Gamma \subset \Theta$
 such that $\Pi_\Gamma\Theta_0 \cap \Pi_\Gamma\Theta_1 = \emptyset$.
\end{thm}
The proofs of Theorems \ref{t23}, \ref{t24} are akin to the proofs of Theorems \ref{t21}, \ref{t22} and are omitted.
\section{\bf Proofs of Theorems of subsection \ref{s2.1} \label{s3}}

{\it Proof of {\sl i.} in Theorem \ref{t3.1}. Preliminary Lemma \ref{l3.1}}. Let $\Lambda_t, t=1,2,$ be the sets
of all measures $P_t, P_t(\Omega_t)<\infty,$ defined
on $\sigma$-field $\frak{B}_t$ of Borel sets in Hausdorff space $\Omega_t$ respectively.
\begin{lemma}\label{l3.1} Let $\Psi_1 \subset \Lambda_1$ and $\Psi_2 \subset \Lambda_2$ be compact in $\tau$-topology.
Let $P_1$ and $P_2$  be cluster points of $\Psi_1$ and $\Psi_2$ respectively. Then the probability measure $P^{(2)} = P_1\otimes P_2$ is
cluster point of $\Psi^{(2)} = \Psi_1\otimes\Psi_2$  in the $\tau$-topology.
\end{lemma}
Lemma \ref{l3.1} is known specialists in measure theory (see 8.10.116 in Bogachev \cite{bo}).
However we do not know publications containing the proof.
\begin{corollary}\label{c1} Let $\Psi\subset \Lambda_1$ be compact in $\tau$-topology. Let $P \in \Lambda_1$ be cluster point of $\Psi$.
Then $P^{(n)} = \otimes_1^n P$ is cluster point of $\Psi^{(n)}= \{Q^{(n)}: Q^{(n)} = \otimes_1^n Q, Q \in \Psi\}$.
\end{corollary}

In Lemma \ref{l3.1} the sets $\Lambda_t$ of all measures can be replaced with the set $\Theta_t$ of all probability measures. Such a statement of Lemma \ref{l3.1} is necessary for the proof of Theorem \ref{t4.1}.

{\it Proof of Lemma \ref{l3.1}}.
By Theorem 2.6 in Ganssler \cite{ga}, a set  $\Psi$ is compact in $\tau$-topology
iff $\Psi$ is sequentially compact in this topology.

Let sequences $P_{tk}\in\Psi_t, t=1,2,$ converge to   $P_t, t=1,2,$ in $\tau$-topology.

Denote $P_k^{(2)} = P_{1k}\otimes P_{2k}$.

Denote $\frak{B}$ the $\sigma$-field in $\Omega=\Omega_1\times\Omega_2$ generated by the product of $\sigma$-fields
$\frak{B}_1$ and $\frak{B}_2$.

It suffices to show that, for any $B \in \frak{B}$ we have
\begin{equation}\label{l1}
\lim_{k\to\infty} P^{(2)}_k(B) = P^{(2)}(B).
\end{equation}
Denote $\bar P^{(2)}_k = P_{1k}\otimes P_1$.

We have
\begin{equation}\label{pus1}
\begin{split}&
|P^{(2)}_k(B) - P^{(2)}(B)| \le |P^{(2)}_k(B) - \bar P^{(2)}_k(B)|\\& + |\bar P^{(2)}_k(B) - P^{(2)}(B)| \doteq |I_{1k}| + |I_{2k}|.
\end{split}
\end{equation}
For all $x \in \Omega_1$ denote $B_x = \{y: (x,y) \in B\}.$

By Fubini Theorem, the sets $B_x, x \in \Omega_2,$ are measurable and
\begin{equation*}
I_{2k} = \int_{\Omega_1} (P_{2k}(B_x) - P_2(B_x)) \,dP_1.
\end{equation*}
Since $ P_{2k}(B_x) \to P_2(B_x)$ as $k \to \infty$ for each $x \in \Omega_1$ we get
\begin{equation}\label{l2}
\lim_{k\to \infty} I_{2k} =0.
\end{equation}
We have
\begin{equation}\label{u2}
I_{1k} = \int_\Omega(P_{2k}(B_x) - P_2(B_x)) \, dP_{1k}.
\end{equation}
By Theorem 2.6 in Ganssler \cite{ga}, there exists probability measure $\mu$ such that $P_{1k}<<\mu$, $P_1<<\mu$ and the functions
$f_k(x) = \frac{dP_{1k}}{d\mu}(x), f(x) = \frac{dP_1}{d\mu}, k = 1,2,\ldots$ are uniformly integrable. By Theorem T22 in Meyer \cite{me},
this implies
\begin{equation}\label{u3}
\lim_{C\to\infty} \sup_k P_k(x: f_k(x) >C) = 0
\end{equation}
Hence, we have
\begin{equation}\label{u4}
 I_{1k} = \int_{\Omega_1}(P_{2k}(B_x) - P_2(B_x)) f_k(x)\, d\mu.
\end{equation}
Denote
\begin{equation*}
f_{1Ck}(x) = f_k(x) 1_{\{f_k(x) > C\}}(x), \quad x \in \Omega_1
\end{equation*}
and $f_{2Ck} = f_k - f_{1Ck}$.

Then
\begin{equation}\label{u6}
   I_{1k} = I_{1k1C} + I_{1k2C}
\end{equation}
where
\begin{equation*}
 I_{1kiC} =  \int_{\Omega_1}(P_{2k}(B_x) - P_2(B_x)) f_{iCk}(x)\, d\mu
\end{equation*}
for $i =1,2$.

Since $ P_{2k}(B_x) \to P_2(B_x)$ as $k \to \infty$ for each $x \in \Omega$, we get
\begin{equation}\label{u8}
I_{1k1C} = o(1) \quad \mbox{as}\quad k\to\infty.
\end{equation}
By (\ref{u3}), we get
\begin{equation}\label{u9}
|I_{1k2C}| \le C\int_{\Omega_1} f_{2Ck}(x) d\mu = o(1) \quad \mbox{as}\quad k \to \infty.
\end{equation}
By (\ref{pus1},\ref{l2},\ref{u6}-\ref{u9}), we get (\ref{l1}).
 This completes the proof of Lemma \ref{l3.1}.

{\it Proof of {\sl i} in Theorem} \ref{t3.1}. Suppose the contrary.
Let $P$ be common cluster point of $\Theta_0$ and $\Theta_1$. Then there exist sequences
$ P_k \in \Theta_0$ and $Q_k \in \Theta_1$ converging to $P \in \Theta$.

For any test $K_n=K_n(X_1,\ldots,X_n)$ we have
\begin{equation}\label{d4}
\lim_{k\to\infty} E_{P_k}[K_n] = E_P [K_n]
\end{equation}
and
\begin{equation}\label{d5}
\lim_{k\to\infty} E_{Q_k}[1-K_n] = E_P [1-K_n].
\end{equation}
By (\ref{d4}) and (\ref{d5}), we get
\begin{equation*}
\alpha(K_n) + \beta(K_n) =1.
\end{equation*}
and come to contradiction. The proof of (\ref{d4}) and (\ref{d5}) is based on approximation $K_n$ simple functions $L_n$ satisfying
(\ref{lu1}) and is omitted. This completes the proof of {\sl i.} in Theorem \ref{l3.1}.

{\it Proofs of {\sl ii} in Theorem \ref{t3.1} and Theorem \ref{t3.3}}. Let $\Theta_0$ be relatively compact and let $\frak{cl}_\tau(\Theta_0)\cap\frak{cl}_\tau(\Theta_1)=\emptyset$.  For any $P \in \Theta$, any measurable sets $A_1,\ldots,A_t \in \frak{B}$ and $\delta_1,\ldots,\delta_t$ denote
 $$
U(P,A_1,\ldots,A_l,\delta_1,\ldots,\delta_l) = \{Q : |P(A_{s}) - Q(A_{s})| < \delta_s, Q \in \Theta, 1 \le s \le t\}.
 $$
 There exists finite covering of $\Theta_0$ by open sets
$
U_i= U(P_i,A_{i1},\ldots,A_{il_i},\delta_{i1},\ldots,\delta_{il_i}), P_i \in \Theta_0
$
such that  $U_i \cap \mbox{cl}_\tau(\Theta_1) = \emptyset, 1 \le i \le l$. The sets $A_{is}, 1\le i \le l, 1 \le s \le l_i$ can be taken as the sets $B_1,\ldots,B_k$ in Theorem \ref{t3.3}.

{\it Proof of} (\ref{m1}). By Theorem \ref{t3.3}, it suffices to verify this statement for the sets of hypotheses $\Theta_{0k}$ and sets of alternatives $\Theta_{1k}$. In the reasoning we can suppose that $B_i\cap B_j=\emptyset, 1 \le i\ne j \le k$.

Denote $S_\delta(x), \delta>0$ the cube in $R^k$  having the length of the side $\delta/2$ and the center $x$.

Denote $\Theta_{0k}(\delta) =  \cup_{x\in\Theta_{0k}} S_{\delta}(x)$ and
$\Theta_{1k}(\delta) =  \cup_{x\in\Theta_{1k}} S_{\delta}(x)$.

There exists $\delta>0$ such that $\Theta_{0k}(2\delta)\cap\Theta_{1k}(2\delta)=\emptyset$.

There exists finite covering $S_\delta(x_{1}),\ldots, S_\delta(x_{l})$ of the set
$\Theta_0$ with $x_{i} \in \Theta_{0k}, 1 \le i \le l$.

Denote
\begin{equation*}
z_{nj} = \sharp\{X_s: X_s \in B_j, 1 \le s \le n\}/n, \quad 1\le j \le k.
\end{equation*}
Define the tests $K_{n\delta}$
\begin{equation*}
\inf_{1 \le i \le l_0}\max_j |z_{nj} - x_{ij}| > 2\delta
\end{equation*}
For each $P \in \Theta_0$ there exists $x_{i}=(x_{i1},\ldots,x_{ik}), 1\le i \le l$ such that $|P(B_j) - x_{ij}| <\delta$.

For any $P \in \Theta_0$ we have
\begin{equation}\label{d6}
\begin{split}&
P(\inf_{1 \le i \le l}\max_{1\le j\le k} |z_{nj} - x_{ij}| > 2\delta) \\& \le  \inf_{1 \le i \le l} P(\max_{1\le j \le k} |z_{nj} - E[z_{nj}] - (x_{ij}- E[z_{nj}])| > 2\delta)\\& \le
P(\max_{1\le j \le k} |z_{nj} - E[z_{nj}]| > \delta)
\le \exp\{-cn\}
\end{split}
\end{equation}
where the last  inequality follows from estimates in the proof of (3.1) in Lemma 3.1 in Groeneboom,  Oosterhoff and  Ruymgaart \cite{gor} with  the constant $c$ depending only on $\delta$ and $k$. As mentioned, the last inequality in (\ref{d6}) is wellknown (see  Barron \cite{ba},
Ermakov \cite{er93} and references therein).

In the case of alternative $P \in \Theta_1$ we have
$$
\inf_{1\le i\le l}\max_{1\le j \le k} |P(B_j) - x_{ij}| > 4\delta.
$$
Therefore
\begin{equation*}
\begin{split}&
\beta(K_{n\delta},P) = P(\inf_{1 \le i \le l}\max_{1\le j\le k} |z_{nj} - x_{ij}| \le 2\delta)\\&\le
l \max_{1 \le i \le l}  P(\max_{1\le j\le k}|z_{nj} - E[z_{nj}] - (x_{ij}- E[z_{nj}])|  \le 2\delta)\\&
\le l P(\max_{1\le j \le k} |z_{nj} - E[z_{nj}]| > 2\delta)
\le l\exp\{-cn\}.
\end{split}
\end{equation*}
 This completes the proof of (\ref{m1}) and Theorem \ref{t3.1}.

{\it Proofs of {\sl i} in Theorem \ref{t3.5}}. The probability measures $P \in \Theta$ are regular. Therefore weak topology separates measures $P \in \Theta$. Hence, by Lemma 2.3 in \cite{ga}, weak topology and $\tau$-topology coincides on $\Theta_0\cup\Theta_1$. It remains to implement Theorem \ref{t3.1}.

{\it Proof of {\sl ii} in Theorem \ref{t3.5}  and Theorem \ref{t3.7}}. Let $\Theta_0$ be relatively compact in weak topology.
There is a finite covering  $\Theta_0$ by  open sets
$$
U_i=U_i(f_{i1},\ldots,f_{im_i},\delta_{i1},\ldots,\delta_{im_i}) =
\left\{  \int_\Omega f_{ij}d(G-P_i) < \delta_{ij}, G \in \Theta, 1\le j \le m_i\right\}
$$
with $f_{ij} \in \bold{C}$, $\delta_{ij} > 0$, $P_i\in \Theta_0$, $1 \le i \le d$, $1 \le j \le m_i$, such that
$\frak{cl}_c(\Theta_1) \cap \cup_{i=1}^m U_i = \emptyset$.

Denote $m = m_1+\ldots+m_d$.

Define the sets
$$
\Psi_t =\left\{v: v = (v_1,\ldots,v_m), v_{m_1+\ldots+m_i+j} = \int_\Omega f_{ij} \, dP, P \in \Theta_t, 1\le i \le d, 1\le j \le m_i\right\}
$$
for $t=0,1$.

For any set $W \subset R^m$ denote $W_{\delta} = \{ u : u = w + v, w \in W, v\in S_\delta(0)\}$
 where $S_\delta(0)$  is the cube having the length of the side
$2\delta$ and the center $0$.

There exists $\delta > 0$ such that $\Psi_{0\delta}\cap \Psi_{1\delta} = \emptyset$.

For any $\epsilon>0, 2\epsilon <\delta$, for each $1 \le i \le m$, $1 \le j \le m_i$ define simple function
$$g_{ij}(x) =\epsilon \sum_{t=1}^l  1_{A_{ijt}}(x) - c_{ij}, \quad
 \,\,x  \in \Omega
 $$
where $A_{ijt} = f_{ij}^{(-1)}((l-t)\epsilon,\infty))$ with $l = l_{ij}= [(d_{ij}-c_{ij})/\epsilon]+1,
d_{ij} = \sup_{x\in\Omega} f_{ij}(x),
c_{ij} = \inf_{x\in\Omega} f_{ij}(x)$. Here $f_{ij}^{(-1)}$ is the inverse function of $f_{ij}$.

It is easy to see that
$$
0 \le g_{ij}(x) - f_{ij}(x) < 2\epsilon \quad \mbox{for all} \quad x \in \Omega.
$$
Define the sets
$$
\Upsilon_t = \left\{v: v = (v_1,\ldots,v_m), v_{m_1+\ldots+m_i+j} = \int_\Omega g_{ij} \, dP,
P \in \Theta_t, 1\le i \le d, 1\le j \le m_i\right\}
$$
for $t=1,2$.

Fix $\delta_1 = \delta - 2\epsilon > 0$. Then
$$
\Upsilon_{0\delta_1}\cap\Upsilon_{1\delta_1} = \emptyset.
$$
Therefore we can take the sets $A_{ijt}, 1 \le i \le m, 1 \le j \le m_i, 1 \le t \le l_{ij}$
 as the sets $B_1,\ldots,B_k$ in Theorem \ref{t3.7}. The remaining part of the reasoning are akin
 to the proofs of Theorems \ref{t3.1}, \ref{t3.3}
 and is omitted.

{\it Proof of Theorem \ref{t3.9}.} In Theorem 3.1 in \cite{ho}, (\ref{m10}) was proved under the certain conditions. These conditions   are satisfied for the sequence of tests $K_{n\delta}$  defined in the proof of (\ref{m1}). It remains to define similar tests $K_{n\delta}$
for the sets $B_1,\ldots,B_k \subset \Omega^m$  defined in Remark 2.1. Note that we also can implement tests $K_{n\delta}$
to  the sets $B_1,\ldots,B_k \subset \Omega^m$   in Remark 2.1 to prove that (\ref{m1}) holds if hypothesis and alternative are not indistinguishable (see Le Cam \cite{le73} and Schwartz \cite{sch}).

{\it Proof of Theorem \ref{t3.10}.}
If $\Omega$ is Polish space, then   $\Theta$ is Polish space as well. Therefore we can define in $\Theta$ the Levi-Prokhorov metric $\rho(P,Q), P,Q \in \Theta$. For each $\delta>0$ define the sets of alternatives
 \begin{equation*}
 \Theta_\delta = \{Q: \inf_{P\in \Theta_0}\rho(Q,P) > \delta, Q \in \Theta\}.
 \end{equation*}
 For any sequence $\delta_n$ define a sequence of tests $L_n = \max_{1\le i \le n} K_{n\delta_i}$ where the tests $K_{n\delta}$ were defined in the proof of Theorem \ref{t3.1}. Applying Theorem \ref{t3.9}, we easily get that there exists  such a sequence $\delta_n\to 0$ as $n \to \infty$ that the sequence of tests $L_n$ is uniformly asymptotically discernible. This completes  the proof of Theorem \ref{t3.10}.
\section{Proofs of Theorems of subsection \ref{s2.3} \label{s4}}
{\it Proof {\sl i}. in Theorem \ref{t4.1}}. Denote $N_n$ the number of atoms $\delta_{X_1},\ldots,\delta_{X_{N_n}}$  of Poisson random measures $\kappa_1,\ldots,\kappa_n$.
Suppose that the sets $\Theta_0$ and $\Theta_1$ have common cluster point $P$ and are not indistinguishable.  There exist sequences
$ P_k \in \Theta_0$ and $Q_k \in \Theta_1$ converging to $P \in \Theta$.

For any test $K_n$, for any $l$, we have
\begin{equation*}
\lim_{k\to\infty} E_{P_k}(K_n|N_n=l) = E_P(K_n|N_n=l),
\end{equation*}
\begin{equation*}
\lim_{k\to\infty} E_{Q_k}(1-K_n|N_n=l) = E_P(1-K_n|N_n=l)
\end{equation*}
and
\begin{equation*}
\lim_{k\to\infty}P_k(N_n = l) =   \lim_{k\to\infty} Q_k(N_n = l) = P(N_n =l).
\end{equation*}
Hence we come to contradiction.

{\it Proof of {\sl ii} in Theorem \ref{t4.1}} is  akin to the proof of {\sl ii} in Theorem \ref{t3.1}. We point out only the differences in the proof of (\ref{d6}). We retain the notation of the proof of Theorem \ref{t3.1} unless otherwise stated.

If $\mu \in \Theta_0$, we have
\begin{equation}\label{d7}
\begin{split}&
P(\inf_{1 \le i \le l_0}\max_{1\le j\le k} |z_{nj} - x_{ij}| > 2\delta)\\& \le \inf_{1 \le i \le l} P(\max_{1\le j \le k} |z_{nj} - \mu(B_j)| > \delta) = \prod_{j=1}^k P(|z_{nj} - \mu(B_j)| > \delta).
\end{split}
\end{equation}
Applying (6.2) in Arcones \cite{ar}, we get
\begin{equation}\label{d8}
\begin{split}&
P(|z_{nj} - \mu(B_j)| > \delta)\le \exp\{-n((\mu(B_j)+\delta)(\log(1+\delta/\mu(B_j)) -1) +\mu(B_j))\}\\& + \exp\{n(\mu(B_j)-\delta)(\log(1-\delta/\mu(B_j)) -1) -\mu(B_j))\}.
\end{split}
\end{equation}
Here we suppose $\log(a) = -\infty$ if $a \le 0$.

This completes the proof of Theorem \ref{t4.1}.

{\it Proof of Theorem \ref{t4.3}}. We begin the reasoning with two Lemmas.

For any $\lambda_1, \lambda_2 \in L_2(\nu)$ define the inner product
$$
(\lambda_1,\lambda_2) = \int_\Omega \lambda_1 \lambda_2\, d\nu
$$
and let $||\lambda_1||^2 = (\lambda_1,\lambda_1)$. For any pair of subspaces $\Gamma_1, \Gamma_2 \subset L_2(\nu)$ denote $\Gamma_1 \oplus \Gamma_2 = \{\lambda: \lambda = \lambda_1 + \lambda_2 , \lambda_1 \in \Gamma_1, \lambda_2 \in \Gamma_2\}$.
Denote  $P_{c}$ the probability measure of Poisson process with the intensity function $\lambda_c(t) \equiv 1, t\in \Omega$.
 \begin{lemma}\label{ik} (see \cite{ik}, Lemma 4) Let $\kappa_1(t)$ and $\kappa_2(t)$ be two Poisson random measures of intensities $\lambda_1, \lambda_2 \in L_2(\nu)$ respectively. Then
 \begin{equation}\label{p1}
 \frac{dP_{\lambda_1}}{dP_c} = \exp\left\{-\int_\Omega(\lambda_1(t) -1) \, d\nu + \int_\Omega \log \lambda_1(t) \, d\kappa_c(t)\right\}
 \end{equation}
and
\begin{equation}\label{p2}
E\left(\frac{dP_{\lambda_1}}{dP_c} \frac{dP_{\lambda_2}}{dP_c}\right) = \exp\left\{\int_\Omega (\lambda_1(t) - 1)(\lambda_2(t) -1)\, d\nu\right\}.
\end{equation}
\end{lemma}
In Ingster and Kutoyants \cite{ik}, (\ref{p2}) was deduced from Campbell formula.
In \cite{ik} the functions $\lambda_1, \lambda_2$ were supposed bounded.
The analysis of the proof of (\ref{p2}) in \cite{ik} shows that it suffices to suppose only that $\lambda_1,\lambda_2 \in L_2(\nu)$.

In the case of the problem of signal detection in Gaussian white noise the reasoning in the proof of Theorem \ref{t1} are the same.
The main difference is that (\ref{p2}) is replaced with
\begin{equation*}
E\left(\frac{dP_{S_1}}{dP_{S_0}} \frac{dP_{S_2}}{dP_{S_0}}\right) =
\exp\left\{\epsilon^{-2}\int_0^1 (S_1(t) - S_0(t))(S_2(t) -S_0(t))\, dt\right\}
\end{equation*}
where
$P_{S_i}, i=0,1,2$ is probability measure of random process $Y_{i\epsilon}$ defined by the stochastic differential equation
$$
dY_{i\epsilon}(t) = S_i(t) dt + \epsilon dw(t).
$$
Denote $\kappa = \kappa_1+\ldots+\kappa_n$.
 \begin{lemma}\label{pu1} Let $\Theta_0  = \{\lambda_0\}, \lambda_0 \in \Theta$ and let $\Theta_1\subset\Theta$. Assume the hypothesis and set of alternatives are distinguishable. Then there exists a finite dimensional subspace $\Gamma \subset\Theta$
such that $\Pi_\Gamma \lambda_0 \notin \Pi_\Gamma \Theta_1$.
\end{lemma}
{\it Proof of Lemma \ref{pu1}}. Suppose the contrary. Suppose for definiteness that $\int_\Omega \lambda_0(t) \, d\nu =1$ and $P_c(\Omega) =1$. Then for any
 sequence $\rho_k, \rho_k \to 0$ as $k \to \infty$  there exists a sequence $\lambda_k \in \Theta_1$ such that, for  all $0\le j < k$, there holds
\begin{equation}\label{p3}
|(\lambda_j -1,\lambda_k - 1)| < \rho_k.
\end{equation}
Define a sequence of Bayes a priori measures $\mu_m$ such that $\mu_m(\lambda_j) = 1/m, 1\le j \le m$. Denote $P_m$ the Bayes a posteriori measure.

The  likelihood ratio for Bayes alternative equals
\begin{equation}\label{p4}
I_m = dP_m/dP_c = \frac{1}{m} \sum_{j=1}^m  \exp\left\{-n\int_\Omega(\lambda_j(t) - \lambda_0(t))\, d\nu + \int_\Omega(\log \lambda_j(t) - 1) d\kappa\right\}.
\end{equation}
By (\ref{p2}) and (\ref{p3}), for the proper choice of $\rho_k, \rho_k \to 0$ as $k \to \infty$, we get
\begin{equation}\label{p5}
\begin{split} &
\mbox{E} [(I_m-dP_{0}/dP_c)^2] = \frac{1}{m^2} \sum_{1\le j_1,j_2\le m}\exp\{n(\lambda_{j_1} - 1,\lambda_{j_2} - 1)\}\\& - \frac{2}{m}\sum_{ j=1}^m\exp\{n(\lambda_{j} - 1,\lambda_0 - 1)\}+1 = o(1)
\end{split}
\end{equation}
as $m \to \infty$.

For fixed $n$, we have
\begin{equation}\label{p6}
\lim_{m \to \infty} \mbox{E} [(I_m-dP_{0}/dP_c)^2] \le \lim_{m\to\infty} (\mbox{Var} [I_m])^{1/2} = 0
\end{equation}
This implies Lemma \ref{pu1}.

{\it Proof of Theorem \ref{t4.3}}. It suffices to prove {\sl i.} Suppose the contrary. For any $\tau \in \Theta_0$ and $\eta \in \Theta_1$ denote $\Gamma_\tau$ and $\Gamma_\eta$ finite dimensional subspaces such that
$\Pi_{\Gamma_\tau}\tau \notin \Pi_{\Gamma_\tau}\Theta_1$  and $\Pi_{\Gamma_\eta}\eta \notin \Pi_{\Gamma_\eta}\Theta_0$. Denote $\Gamma_{\tau\eta} = \Gamma_\tau\oplus\Gamma_\eta$. Denote $\Upsilon_\tau$ and $\Upsilon_\eta$ the  linear spaces generated the vectors $\tau$ and $\eta$ respectively.

For any $\tau \in \Theta_0$ and $\eta \in \Theta_1$ there exists $\gamma \in \Pi_{\Gamma_{\tau\eta}}\Theta_0 \cap \Pi_{\Gamma_{\tau\eta}}\Theta_1$.

Let us show that there exist sequences of points $\tau_i \in \Theta_0, \eta_i \in \Theta_1$ and   finite dimensional  subspaces $\Gamma_i$ such that

{\it i.}  $\lambda_c \equiv 1 \in \Gamma_0$.

{\it ii.} $\Gamma_i = \Gamma_{i-1}  \oplus \Gamma_{\tau_{i}\eta_{i}}\oplus\Upsilon_{\tau_{i-1}}\oplus\Upsilon_{\eta_{i-1}}$.

{\it iii.} there  exists a sequence of points $\gamma_{ii} = \Pi_{\Gamma_i}\tau_{i+1}=\Pi_{\Gamma_i}\eta_{i+1}$.

{\it iv.} for each $i$ there exists $z_i \in  \Gamma_i$ such that $\gamma_{ij} = \Pi_{\Gamma_i} \gamma_{jj} \to z_i$ as $j \to \infty$.

{\it v.} $$ \int_\Omega \tau_i dP_c \to C,\quad \int_\Omega \eta_i dP_c \to C\quad \mbox{\rm as} \quad i \to \infty $$

Here $ \Upsilon_{\tau_0} = \Upsilon_{\eta_0}=\emptyset$.

Denote $\Gamma_{ci} = \Gamma_{\tau_{i}} \oplus \Gamma_{\eta_{i}}\oplus\Upsilon_{\tau_{i-1}}\oplus\Upsilon_{\eta_{i-1}}.$

We can define sequences $\tau_i \in \Theta_0$ and $\eta_i \in \Theta_1$ satisfying {\it i.-iii.} by induction.

Let $\tau_1 \in \Theta_0$ and $\eta_1 \in \Theta_1$. Denote $\Gamma_1 =\Upsilon_\lambda \oplus \Gamma_{\tau_1\eta_1}$.
Let $\gamma_{11} \in \Pi_{\Gamma_1}\Theta_0\cap \Pi_{\Gamma_1}\Theta_1$. Define $\tau_2\in \Theta_0$ and $\eta_2\in \Theta_1$ as arbitrary points such that
$\gamma_{11}=\Pi_{\Gamma_1}\tau_2 = \Pi_{\Gamma_1}\eta_2$ and so on.

Using these sequences $\tau_i, \eta_i$ satisfying {\it i.-iii.} we find a subsequence satisfying {\it i.-iv.} on
the base of the following procedure.
We choose a subsequence $\gamma_{i_{s_1}i_{s_1}}$ such that $\gamma_{i_{s_1}1}$ converges to some point $z_1 \in \Gamma_1$. After that we choose from these subsequence a subsequence $\gamma_{i_{s_2}i_{s_2}}$ such that $\gamma_{i_{s_2}2}$ converges to some point $z_2 \in \Gamma_{2_{s_1}}$ and so on. The sequences of points $\tau_{i_{s_i}}$,  $\eta_{i_{s_i}}$, $\gamma_{i_{s_i},i_{s_i}}$ and subspaces $\Gamma_{i_{s_{i-1}}}$  satisfy {\it i- iv}.

The same procedure we can make such that {\it v.} holds.

Let sequences $\tau_{i}$,  $\eta_{i}$, $\gamma_{ii}$ and subspaces $\Gamma_{i}$ will satisfy {\it i- v}. Define the sets $ \Theta_{0m} = \{\tau_1,\ldots,\tau_m\} \subset \Theta_0$ and $ \Theta_{1m} = \{\eta_1,\ldots,\eta_m\} \subset \Theta_1$. For testing the hypothesis $\bar H_{0m}: \mu \in \Theta_{0m}$ versus $\bar H_{1m}: \mu \in \Theta_{1m}$ we introduce Bayes a priori measures $\mu_{0m}$ and $\mu_{1m}$ such that $\mu_{0m}(\tau_i) = m^{-1}$ and $\mu_{1m}(\eta_i) = m^{-1}$ with $1 \le i \le m$.

By an appropriate choice of subsequence $i_k$ we can make the differences $\gamma_{ij}- z_{ii},    \int_\Omega \tau_i dP_c - C$ and $\int_\Omega \eta_i dP_c- C$ negligible for further  estimates. Thus we shall assume that $\gamma_{ij} = z_i, 1 \le i < \infty, i \le j < \infty$ in the further reasoning. This allows us to choose a system of coordinates such that
\begin{equation*} \begin{split}& z_{i} = \sum_{s=0}^i a_s\psi_s,\\&
\tau_{i+1} = z_{i} + a_{1i}\psi_{i+1}+b_{1i}\zeta_{1i} \quad\mbox{\rm and} \quad \eta_{i+1} = z_{i} + a_{2i}\psi_{i+1}+b_{2i}\zeta_{1i} + c_{2i}\zeta_{2i}
\end{split}
\end{equation*}
where $\psi_i, \zeta_{1i}, \zeta_{2i}, 1\le i < \infty$ are orthogonal functions and $\psi_0\equiv 1$.

Denote $\nu_{lm}, l=0,1$ Bayes a posteriori probability measure having a priori measure $\mu_{lm}$. Denote $\pi_{lm} = d\nu_{lm}/dP_0$.

We have
\begin{equation}\label{p7}
\pi_{0m} - \pi_{1m} = \frac{1}{m}\sum_{i=1}^m J_i
\end{equation}
where
\begin{equation*}
\begin{split}&
J_i = \exp\left\{-n\sum_{s=1}^i a_s\int_\Omega \psi_s(t)\,d\nu  + na_{1i}\int_\Omega \psi_{i+1}(t)\,d\nu\right.\\&\left. + nb_{1i} \int_\Omega \zeta_{1i}(t)\,d\nu  + \int_\Omega \log \tau_{i+1}(t)\, d\kappa\right\}\\& -
\exp\left\{-n\sum_{s=1}^i a_s\int_\Omega \psi_s(t)\,d\nu  + na_{2i}\int_\Omega \psi_{i+1}(t)\,d\nu + nb_{2i} \int_\Omega \zeta_{1i}(t)\,d\nu\right.\\& \left. + nc_{2i} \int_\Omega \zeta_{2i}(t)\,d\nu + \int_\Omega \log \eta_{i+1}(t)\, d\kappa\right\}.
\end{split}
\end{equation*}
Let $k_1 < k_2$. Then, by (\ref{p2}), we get
\begin{equation}\label{p8}
\begin{split}&
E[J_{i_1}J_{i_2}] = \exp\{n(\tau_{i_1},\tau_{i_2})\} - \exp\{n(\tau_{i_1},\eta_{i_2})\}- \exp\{n(\eta_{i_1},\tau_{i_2})\} + \exp\{n(\eta_{i_1},\eta_{i_2})\}\\&=
\exp\left\{\sum_{s=1}^{i_1} na_s^2\right\}
(\exp\{na_{1i_1}a_{i_1}\} - \exp\{na_{1i_1}a_{i_1}\}\\& -\exp\{na_{2i_1}a_{i_1}\} +\exp\{na_{2i_1}a_{i_1}\})  = 0.
\end{split}
\end{equation}
We also have
\begin{equation}\label{p9}
EJ_i^2 \le 2 \exp\left\{\sum_{s=1}^i na_s^2\right\}(\exp\{na_{1i}^2 + nb_{1i}^2\} + \exp\{na_{2i}^2 + nb_{2i}^2 + nc_{2i}^2\}).
\end{equation}
By (\ref{p8}) and (\ref{p9}), applying the Cauchy inequality, we get
\begin{equation*}
\lim_{m \to \infty} (E|\pi_{1m} - \pi_{2m}|)^2 \le  \lim_{m \to \infty} E(\pi_{1m} - \pi_{2m})^2 = 0.
\end{equation*}
Since $m$ does not depend on $n$ we come to contradiction. This completes the proof of Theorem \ref{t4.3}.

\vskip 0.3cm
 Mechanical Engineering Problems Institute\\
Russian Academy of Sciences\\
Bolshoy pr.,V.O., 61\\
St.Petersburg\\
Russia\\
e-mail: erm2512@mail.ru

\end{document}